\documentclass[11pt,a4paper]{article}

\usepackage{amscd,amssymb,amsmath}
\usepackage{graphicx}

\usepackage[latin1]{inputenc}
\usepackage{color}

\date{}

\newtheorem{theorem}{Theorem}[section]

\newtheorem{definition}{Definition}[section]

\newtheorem{lemma}[theorem]{Lemma}

\def\pf{\noindent {\it Proof.} }

\def\qed{\hfill \nopagebreak\rule{5pt}{8pt}}
\def\pf{\noindent {\it Proof.} }

\title{Matching Energy of Unicyclic and Bicyclic
Graphs with a Given Diameter}

\begin{document}

\maketitle

\noindent Lin Chen, Jinfeng Liu, Yongtang Shi\footnote{The
corresponding author.}

\vspace{0,4cm} \noindent  Center for Combinatorics and LPMC-TJKLC,
Nankai University,
Tianjin 300071, P.R. China\\

\noindent E-mails:  chenlin1120120012@126.com, ljinfeng709@163.com,
shi@nankai.edu.cn

\begin{abstract}
Gutman and Wagner proposed the concept of matching energy (ME) and
pointed out that the chemical applications of ME go back to the
1970s. Let $G$ be a simple graph of order $n$ and
$\mu_1,\mu_2,\ldots,\mu_n$ be the roots of its matching polynomial.
The matching energy of $G$ is defined to be the sum of the absolute
values of $\mu_{i}\ (i=1,2,\ldots,n)$. In this paper, we
characterize the graphs with minimal matching energy among all
unicyclic and bicyclic graphs with a given diameter $d$.
\end{abstract}


\section{Introduction}
In this paper, all graphs under our consideration are finite,
connected, undirected and simple. For more notations and terminology
that will be used in the sequel, we refer to \cite{graphBondy2008}.
Let $G$ be a simple undirected graph with order $n$ and $A(G)$ be
the adjacency matrix of $G$.  The characteristic polynomial of $G$,
denoted by $\phi(G)$, is defined as
\begin{equation*}
\phi(G)=\text{det}(xI-A(G))=\sum\limits_{i=0}^n a_i(G)x^{n-i},
\end{equation*}
where $I$ is the identity matrix of order $n$. The roots of the
equation $\phi(G)=0$, denoted by
$\lambda_1,\lambda_2,\ldots,\lambda_n$, are the eigenvalues of
$A(G)$.  The energy of $G$, denoted by $E(G)$, is defined as the sum
of the absolute values of the eigenvalues of $A(G)$, that is,
$$E(G)=\sum_{i=1}^n|\lambda_i|.$$
The concept of the energy of simple undirected graphs was introduced
by Gutman in \cite{G1} and now is well-studied. For more results
about graph energy, we refer the readers to recent papers \cite{DM,
DMG,DLS,LWS, MMZ}, two surveys \cite{G2,G3} and the book \cite{LSG}.
There are various generalizations of graph energy, such as
Randi\'{c} energy \cite{BB0, DS}, Laplacian energy \cite{DGCZ},
distance energy \cite{SMHP}, incidence energy \cite{BB, BG}, energy
of matrices \cite{GFAHG} and energy of a polynomial \cite{MVI}, etc.

Let $G$ be a simple graph with $n$ vertices and $m$ edges. Denote by
$m_k(G)$ the number of $k$-matchings($=$ the number of selections of
$k$ independent edges $=$ the number of $k$-element independent edge
sets) of $G$. Specifically, $m_1(G)=m$ and $m_k(G)=0$ for
$k>\lfloor\frac{n}{2}\rfloor$ or $k<0$. It is both consistent and
convenient to define $m_0(G)=1$. The matching polynomial of the
graph $G$ is defined as
\begin{equation}\label{1}
\alpha(G)=\alpha (G,\mu)=\sum\limits_{k\geq 0} (-1)^k
m_k(G)\mu^{n-2k}.
\end{equation}

Recently, Gutman and Wagner \cite{GW} defined the matching energy of
a graph $G$ based on the zeros of its matching polynomial
\cite{EJF,gutmanmatch1979}.

\begin{definition}
Let $G$ be a simple graph with order $n$, and
$\mu_1,\mu_2,\ldots,\mu_n$ be the zeros of its matching polynomial.
Then,
\begin{equation}\label{2}
ME(G)=\sum\limits_{i=1}^n|\mu_i|.
\end{equation}
\end{definition}
Moreover, Gutman and Wagner \cite{GW} pointed out that the matching
energy is a quantity of relevance for chemical applications. They
arrived at the simple relation: $$TRE(G) = E(G) - ME(G),$$ where
TRE($G$) is the so-called ``topological resonance energy" of $G$.
About the chemical applications of matching energy, for more details
see \cite{gutmanmt1975,aihara1976,gutmanmt1977}.

For the coefficients $a_i(G)$ of $\phi(G)$, let $b_i(G)=|a_i(G)|,
i=0,1,\ldots,n$. Note that $b_0(G)=1$, $b_1(G)=0$, and $b_2(G)$ is
the number of edges of $G$. For convenience, let $b_i(G)=0$ if
$i<0$. In \cite {gutman1977,HY}, we have
\begin{equation}\label{3}
E(G)=\frac{1}{2\pi} \int_{-\infty}^\infty
\frac{dx}{x^2}\ln\Big[\Big(\sum\limits_{j=0}^{\lfloor
\frac{n}{2}\rfloor}
b_{2j}(G)x^{2j}\Big)^2+\Big(\sum\limits_{j=0}^{\lfloor
\frac{n}{2}\rfloor} b_{2j+1}(G)x^{2j+1}\Big)^2\Big].
\end{equation}
Thus $E(G)$ is a monotonically increasing function of $b_i(G),
i=0,1,\ldots,n$.

Being similar to Eq.(\ref{3}), the matching energy also has a
beautiful formula as follows\cite{GW}. Eq.(\ref{4}) could be
considered as the definition of matching energy, in which case
Eq.(\ref{2}) would become a theorem.

\begin{theorem}\label{thm1}
Let $G$ be a simple graph of order $n$, and $m_k(G)$ be the number
of its $k$-matchings, $k = 0, 1,
2,\ldots,\lfloor\frac{n}{2}\rfloor$. The matching energy of $G$ is
given by
\begin{equation}\label{4}
ME=ME(G)=\frac{2}{\pi}\int^{\infty}_{0}\frac{1}{x^2}\ln\Big[
\sum\limits_{k\geq 0}m_k(G)x^{2k}\Big]dx. \qquad\qquad
\end{equation}
\end{theorem}

By Eq.(\ref{4}) and the monotony of the function logarithm, we can
define a \emph{quasi-order} ``$\succeq$" as follows: If two graphs
$G_1$ and $G_2$ have the same order and size, then
$$G_{1}\succeq G_{2}\Longleftrightarrow m_k(G_1)\geq m_k(G_2)\quad \text{for all}\ k.$$
If $G_1\succeq G_2$ and there exists some $k$ such that $m_k(G_1)>
m_k(G_2)$, then we write $G_1\succ G_2$. Clearly, $G_1\succ G_2
\Longrightarrow ME(G_1)> ME(G_2).$

Notice that when $ME(G_1)> ME(G_2)$, we may not deduce that
$G_1\succ G_2$. However, if $G$ is any simple connected graph with
$n$ vertices other than $S_n$, where $S_n$ is a star of order $n$,
then not only $ME(G)>ME(S_n)$ \cite{GW} but also $G\succ S_n$. Based
on the quasi-order, there are some more extremal results on matching
energy of graphs \cite{CS, JLS, JIMA, LY}.

In this paper, we characterize the graphs with minimal matching
energy among all unicyclic and bicyclic graphs with a given diameter
$d$.

\section{Preliminaries}

The following result gives two fundamental identities for the number
of $k$-matchings of a graph (see \cite{EJF, gutmanmatch1979}).
\begin{lemma}\label{lem1}
Let $G$ be a simple graph, $e=uv$ be an edge of $G$, and
$N(u)=\{v_1(=v),v_2,\ldots,v_j\}$ be the set of all neighbors of $u$
in $G$. Then we have
\begin{equation}\label{5}
m_k(G)=m_k(G-uv)+m_{k-1}(G-u-v),
\end{equation}
\begin{equation}\label{6}
m_k(G)=m_k(G-u)+\sum\limits_{i=1}^j m_{k-1}(G-u-v_i).
\end{equation}
\end{lemma}

From Lemma \ref{lem1}, we know that $m_k(P_1\cup G)=m_k(G)$. And we
can also obtain that
\begin{lemma}\label{lem2}
Let $G$ be a simple graph and $H$ be a subgraph(resp. proper
subgraph) of $G$. Then $G\succeq H$(resp. $\succ H$).
\end{lemma}

A connected graph with $n$ vertices and $n$ edges is called a
unicyclic graph. Obviously, a unicyclic graph has exactly one cycle.
A connected graph with $n$ vertices and $n+1$ edges is called a
bicyclic graph. Let $\mathcal {U}(n)$ be the class of connected
unicyclic graphs with $n$ vertices, $\mathcal {U}(n,d)$ be the class
of unicyclic graphs with $n$ vertices and diameter $d$, where $1\leq
d\leq n-2$. Let $\mathcal {B}(n)$ be the class of bicyclic graphs
with $n$ vertices and $\mathcal {B}(n,d)$ be the class of bicyclic
graphs in $\mathcal {B}(n)$ with diameter $d$, where $2\leq d\leq
n-2$. Let $P_n$ be the path with $n$ vertices and $K_n$ be the
complete graph with $n$ vertices.

\begin{figure}[h,t,b,p]
\begin{center}
\includegraphics[scale = 0.7]{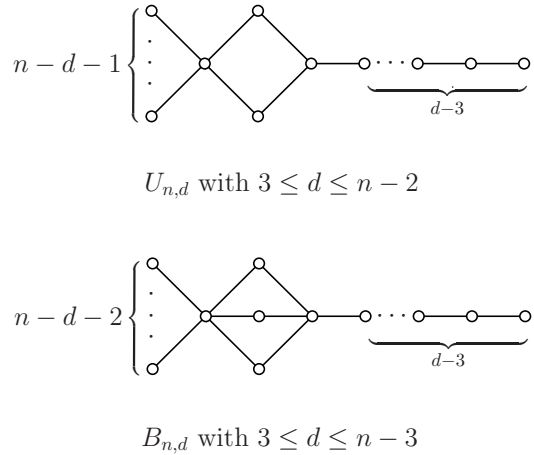}
\caption{The graphs $U_{n,d}$ and $B_{n,d}$.}\label{minimalgraph}
\end{center}
\end{figure}

When $d=1$, $n=3$, $K_3$ is the unique graph in $\mathcal {U}(3,1)$.
When $d=1$, $n\geq 4$, $\mathcal {U}(n,1)$ contains no graphs. When
$d=2$, $n=4$, $\mathcal {U}(4,2)$ has two graphs $G_{4,2}^1$ and
$G_{4,2}^2$ (see Figure \ref{U_4,2}). Clearly, $ME(G_{4,2}^2)>
ME(G_{4,2}^1)$, i.e., $G_{4,2}^1$ is the unique graph with minimal
matching energy in $\mathcal {U}(4,2)$. When $d=2, n\geq 5$, the
graph obtained by attaching $n-3$ pendant vertices to a vertex of a
triangle is the unique graph in $\mathcal {U}(n,2)$. Thus, we just
consider the case in which $3\leq d\leq n-2$. In section $3$ of our
paper, we will prove that for $3\leq d\leq n-2$, the graph $U_{n,d}$
is the unique graph in $\mathcal {U}(n,d)$ with minimal matching
energy, where the graph $U_{n,d}$ is shown in Figure
\ref{minimalgraph}.

\begin{figure}[h,t,b,p]
\begin{center}
\includegraphics[scale = 0.7]{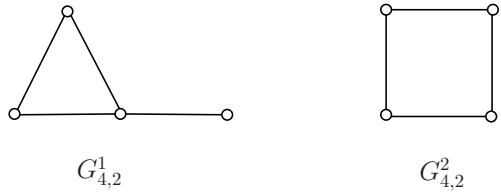}
\caption{The two graphs in $\mathcal {U}(4,2)$.}\label{U_4,2}
\end{center}
\end{figure}

When $d=2$, $n \geq 6$, $\mathcal {B}(n,2)$ has two graphs
$G_{n,2}^1$ and $G_{n,2}^2$(see Figure \ref{B_n,2}). By Lemma
\ref{lem1} and simple calculation, we can get $G_{n,2}^1 \succ
G_{n,2}^2$, hence $G_{n,2}^2$ is the unique graph in $\mathcal
{B}(n,2)$ with minimal matching energy. Therefore, we only consider
the case in which $3\leq d\leq n-2$. In section $4$, we will show
that $B_{n,d}$ is the unique graph with minimal matching energy for
$3\leq d\leq n-3$, where the graph $B_{n,d}$ is shown in Figure
\ref{minimalgraph}. Furthermore, we also pay our attention to the
case $d=n-2$.

\begin{figure}[h,t,b,p]
\begin{center}
\includegraphics[scale = 0.7]{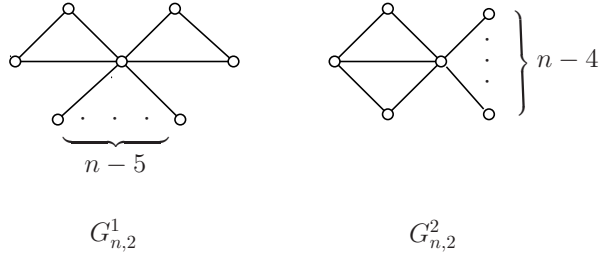}
\caption{The two graphs in $\mathcal {B}(n,2)$.}\label{B_n,2}
\end{center}
\end{figure}

Let $\mathcal {T}(n,d)$ be the class of trees with $n\geq2$ vertices
and diameter $d$, where $1\leq d\leq n-1$. If $T\in \mathcal
{T}(n,1)$, then $T=P_2$. For $1\leq d\leq n-1$, let $T_{n,d}$ denote
the graph obtained by attaching $n-d$ pendent vertices to an end
vertex of $P_d$. Specially, $T_{n,1}=T_{n,2}=S_n$. Obviously,
$T_{n,2}$ is the unique tree in $\mathcal {T}(n,2)$ and
$T_{n,n-1}=P_n$ is the unique tree in $\mathcal {T}(n,n-1)$.

Let $G_1, G_2$ be two graphs with $n$ vertices. Now we introduce a
\emph{quasi-order} $\succeq_1$ defined in \cite{LSG}: If
$b_i(G_1)\geq b_i(G_2)$ for all $i\geq 0$, then we write $G_1
\succeq_1 G_2$. If $G_1 \succeq_1 G_2$ and there exists an $i_0$
such that $b_{i_0}(G_1)>b_{i_0}(G_2)$, then we write $G_1 \succ_1
G_2$. The following lemmas are relevant results on this quasi-order.

\begin{lemma}[\cite{GP,ZL}]\label{PPP}
For $2\leq i\leq \lfloor\frac{n}{2}\rfloor $ and $n\geq 4$,
$$P_n \succ_1 P_i \cup P_{n-i} \succ_1 P_1 \cup P_{n-1}.$$
\end{lemma}

\begin{lemma}[\cite{gutman1977}]\label{PTS}
For $3\leq d\leq n-2$, $P_n\succeq_1 T_{n,d}\succeq_1 S_n$.
\end{lemma}

\begin{lemma}[\cite{YY}]\label{TT}
Let $T\in \mathcal {T}(n,d)$ and $T\neq T_{n,d}$. Then $T\succ_1
T_{n,d}$.
\end{lemma}

\begin{lemma}[\cite{LZ}]\label{TdTd_0}
If $d>d_0\geq 3$, then $T_{n,d}\succ_1 T_{n,d_0}$.
\end{lemma}

\begin{lemma}[\cite{YZ}]\label{TTT}
For $2\leq d_1 \leq n_1-2$, we have $T_{n_1,d_1}\cup T\succeq_1
T_{n_1+n_2-1,d_1+d_2}$, where $T=T_{n_2,d_2}$ if $2\leq d_2 \leq
n_2-2$, and $P_2$ if $n_2=2$ and $d_2=1$.
\end{lemma}

If $G$ is an acyclic graph, then \cite{GP} $b_{2k}(G)=m_k(G)$ and
$b_{2k+1}(G)=0$ for all $k$. Thus, the quasi-order $\succ_1$ (resp.
$\succeq_1$) in Lemmas \ref{PPP}--\ref{TTT} can be replaced by
$\succ$ (resp. $\succeq$), and the results also work.

By Lemma \ref{lem1} and the definition of the quasi-order $\succeq$,
it is easy to see that the following lemma holds.

\begin{lemma}\label{lem9}
Let $G$, $G'\in \mathcal {U}(n)$ and $uv$(resp. $u'v'$) be a pendant
edge with the pendant vertex $u$(resp. $u'$) of the graph $G$(resp.
$G'$). If $G-u\succeq G'-u'$ and $G-u-v\succ G'-u'-v'$, or $G-u\succ
G'-u'$ and $G-u-v\succeq G'-u'-v'$, then $G\succ G'$.
\end{lemma}

The following lemmas will be needed in our paper, which are obtained
based on the previous results.

\begin{lemma}\label{BUT}
For $3\leq d \leq n-2$, $B_{n,d}\succ U_{n,d}\succ T_{n,d}$.
\end{lemma}
\pf Since $U_{n,d}$ is a proper subgraph of $B_{n,d}$, then by Lemma
\ref{lem2}, we can get $B_{n,d}\succ U_{n,d}$. Similarly,  we also
have $U_{n,d}\succ T_{n,d}$. \qed

\begin{lemma}\label{Udd_0}
For $3\leq d_0 < d \leq n-2$, $U_{n,d}\succ U_{n,d_0}$.
\end{lemma}
\pf By Lemmas \ref{lem1}, \ref{lem2} and \ref{TdTd_0},
\begin{eqnarray*}
m_k(U_{n,d}) &\geq& m_k(U_{n-1,d-1})+m_{k-1}(T_{n-2,d-2})\\
             &\geq& m_k(U_{n-1,d-1})+m_{k-1}(T_{d-1,d-3})\\
             &=& m_k(U_{n,d-1}).
\end{eqnarray*}
Furthermore, $m_2(U_{n,d})> m_2(U_{n,d-1})$. It follows that
$U_{n,d}\succ U_{n,d-1}$. Therefore, $U_{n,d}\succ U_{n,d-1}\succ
\cdots \succ U_{n,d_0}$. \qed

Similarly, we have

\begin{lemma}\label{Bdd_0}
For $3\leq d_0 < d \leq n-2$, $B_{n,d}\succ B_{n,d_0}$.
\end{lemma}

\section{Unicyclic graphs with a given diameter}
Now we consider the minimal matching energy of graphs in $\mathcal
{U}(n,d)$ with $3\leq d\leq n-2$. We first discuss the case $d=n-2$.
\begin{lemma}\label{lem13}
Let $G\in \mathcal {U}(n,n-2)$ with $n\geq 8$ and $G\neq U_{n,n-2}$.
Then $G\succ U_{n,n-2}$.
\end{lemma}
\pf We will prove the lemma by induction on $n$.

If $n=8$, then $G$ is isomorphic to one of the following graphs (see
Figure \ref{U_8,6}).

\begin{figure}[h,t,b,p]
\begin{center}
\includegraphics[scale = 0.6]{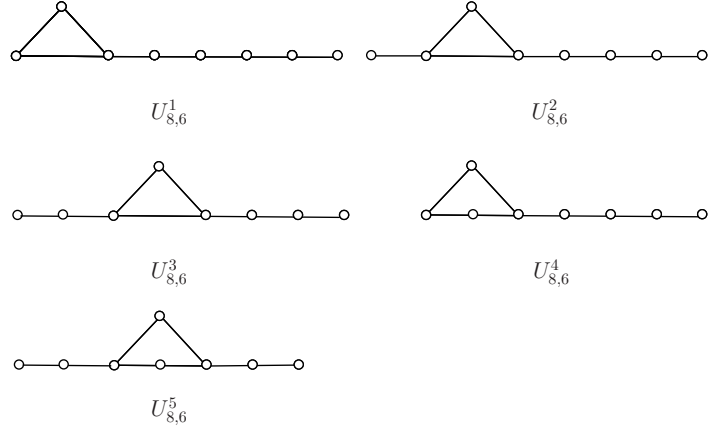}
\caption{The graphs in $\mathcal {U}(8,6)$ except for
$U_{8,6}$.}\label{U_8,6}
\end{center}
\end{figure}

It is easy to get
$$\alpha(U_{8,6}^1)=\mu^8-8\mu^6+19\mu^4-13\mu^2+1;~~~
\alpha(U_{8,6}^2)=\mu^8-8\mu^6+18\mu^4-11\mu^2+1;$$
$$\alpha(U_{8,6}^3)=\mu^8-8\mu^6+18\mu^4-12\mu^2+1;~~~
\alpha(U_{8,6}^4)=\mu^8-8\mu^6+19\mu^4-14\mu^2+2;$$
$$\alpha(U_{8,6}^5)=\mu^8-8\mu^6+18\mu^4-12\mu^2+2;~~~
\alpha(U_{8,6})=\mu^8-8\mu^6+18\mu^4-11\mu^2.~~~~~$$

And then it is obvious that $U_{8,6}^i \succ U_{8,6}$ for
$i=1,2,3,4,5$, i.e., $U_{8,6}$ is the unique graph in $\mathcal
{U}(8,6)$ with minimal matching energy.

If $n=9$, then $G$ is isomorphic to one of the following graphs (see
Figure \ref{U_9,7}).

\begin{figure}[h,t,b,p]
\begin{center}
\includegraphics[scale = 0.7]{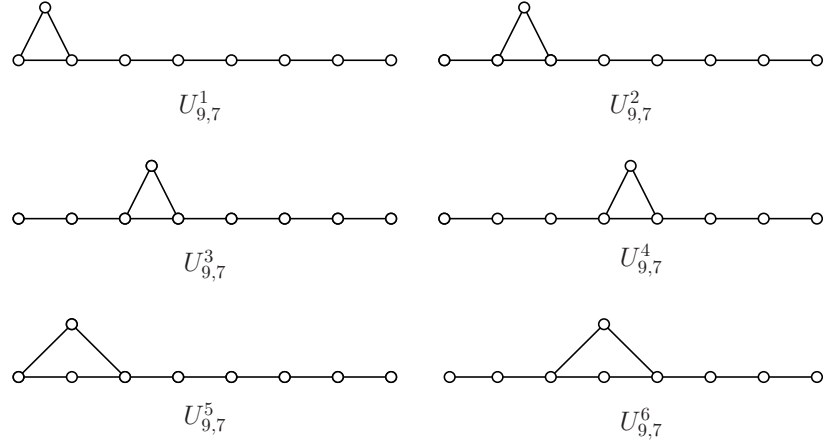}
\caption{The graphs in $\mathcal {U}(9,7)$ except for
$U_{9,7}$.}\label{U_9,7}
\end{center}
\end{figure}

We can obtain that
$$\alpha(U_{9,7}^1)=\mu^9-9\mu^7+26\mu^5-26\mu^3+6\mu;~~~
\alpha(U_{9,7}^2)=\mu^9-9\mu^7+25\mu^5-23\mu^3+5\mu;$$
$$\alpha(U_{9,7}^3)=\mu^9-9\mu^7+25\mu^5-24\mu^3+6\mu;~~~
\alpha(U_{9,7}^4)=\mu^9-9\mu^7+25\mu^5-24\mu^3+5\mu;$$
$$\alpha(U_{9,7}^5)=\mu^9-9\mu^7+26\mu^5-27\mu^3+8\mu;~~~
\alpha(U_{9,7}^6)=\mu^9-9\mu^7+25\mu^5-24\mu^3+6\mu;$$ along with
$$\alpha(U_{9,7})=\mu^9-9\mu^7+25\mu^5-23\mu^3+4\mu.$$

It now immediately follows that $U_{9,7}^i \succ U_{9,7}$ for
$i=1,2,3,4,5,6$, i.e., $U_{9,7}$ is the unique graph in $\mathcal
{U}(9,7)$ with minimal matching energy.

Now suppose that the result holds for graphs in $\mathcal
{U}(n-1,n-3)$ and $\mathcal {U}(n-2,n-4)$. Let  $G\in \mathcal
{U}(n,n-2)$ and $G\neq U_{n,n-2}$, where $n\geq 10$.

Let $u$(resp. $u'$) be a pendant vertex, adjacent to $v$(resp.
$v'$), which has the largest distance to a vertex on the unique
cycle of $G$(resp. $U_{n,n-2}$). Then the degree of $v$ is 2. So is
$v'$. Hence $G-u \in \mathcal {U}(n-1,n-3)$, $G-u-v\in \mathcal
{U}(n-2,n-4)$ and $U_{n,n-2}-u'=U_{n-1,n-3}$,
$U_{n,n-2}-u'-v'=U_{n-2,n-4}$.

Since $G\neq U_{n,n-2}$, we have either $G-u\neq U_{n-1,n-3}$ or
$G-u-v\neq U_{n-2,n-4}$. By the induction hypothesis, we have
$G-u\succ U_{n-1,n-3}$ and $G-u-v\succeq U_{n-2,n-4}$, or
$G-u\succeq U_{n-1,n-3}$ and $G-u-v\succ U_{n-2,n-4}$. By Lemma
\ref{lem9}, $G\succ U_{n,n-2}$. \qed

\begin{theorem}\label{th2}
Let $G\in \mathcal {U}(n,d)$ with $n \geq 8$, $3\leq d\leq n-2$ and
$G\neq U_{n,d}$. Then $ME(G)>ME(U_{n,d})$.
\end{theorem}
\pf We prove the result by induction on $n-d$.

When $n-d=2$, by Lemma \ref{lem13}, we have $G\succ U_{n,d}$. Let
$t\geq 3$ and suppose that the result holds for $n-d<t$. Now suppose
that $n-d=t$. Let $u'$ be the vertex of degree 3 in $U_{n,d}$ and
$v'$ be a vertex on the quadrangle that is adjacent to $u'$. By
Lemma \ref{lem1},
\begin{eqnarray*}
m_k(U_{n,d}) &=& m_k(U_{n,d}-u'v')+m_{k-1}(U_{n,d}-u'-v')\\
             &=& m_k(T_{n,d})+m_{k-1}(P_{d-3}\cup S_{n-d+1}).
\end{eqnarray*}

For $C_n$, $m_k(C_n)=m_k(P_n)+m_{k-1}(P_{n-2})$. By Lemma \ref{PTS},
$P_n \succ T_{n,d}$. And by Lemmas \ref{PPP} and \ref{PTS}, $P_{n-2}
\succeq P_{d-3}\cup P_{n-d+1} \succ P_{d-3}\cup S_{n-d+1}$. Thus
$C_n \succ U_{n,d}$.  Therefore, we may suppose that the unique
cycle of $G$ is $C_r$ with $r<n$. Let $P(G)=v_0v_1\ldots v_d $ be a
diametrical path of $G$. Then one of $v_0$ and $v_d$ must be a
pendant vertex.

{\bf{Case 1}} All pendant vertices are on $P(G)$.

Since $t=n-d\geq 3$, then $|V(P(G))|=d+1\leq n-2$. Thus there are at
least two adjacency vertices, say $u$ and $v$, on $C_r$ which lie
outside $P(G)$ such that $G-uv\in \mathcal {T}(n,d_1)$, and
$G-u-v\in \mathcal {T}(n-2,d_2)$, where $d_1,d_2 \geq d$. By Lemmas
\ref{TT} and \ref{TdTd_0}, $G-uv \succeq T_{n,d_1}\succeq T_{n,d}$,
$G-u-v \succeq T_{n-2,d_2}\succeq T_{n-2,d}$. We also have
$T_{n-2,d}\succ P_{d-3}\cup T_{n-d+1,3} \succ P_{d-3}\cup S_{n-d+1}$
by Lemmas \ref{lem2} and \ref{PTS}. Hence, $G-u-v\succ P_{d-3}\cup
S_{n-d+1}$.

Note that $m_k(G)=m_k(G-uv)+m_{k-1}(G-u-v)$, so $m_k(G)\geq
m_k(U_{n,d})$ for all $k$. Moreover, since $T_{n-2,d}\succ
P_{d-3}\cup S_{n-d+1}$, there exists some $k_0$ such that
$m_{k_0-1}(G-u-v)>m_{k_0-1}(P_{d-3}\cup S_{n-d+1})$, i.e.,
$m_{k_0}(G)>m_{k_0}(U_{n,d})$. Thus $G\succ U_{n,d}$.

{\bf{Case 2}} There is at least one pendant vertex outside $P(G)$.

Let $u'$ be a pendant vertex of $U_{n,d}$ adjacent to the vertex
$v'$ of degree $n-d+1$. Then $U_{n,d}-u'=U_{n-1,d}$, and
$U_{n,d}-u'-v'=(n-d-2)P_1\cup T_{d,d-2}.$

{\bf{Subcase 2.1}} There is a pendant vertex $u$ outside $P(G)$ such
that its neighbor $v$ lies on $C_r$.

 Since $u$ outside $P(G)$, then $G-u\in \mathcal {U}(n-1,d)$. Consequently, by the induction hypothesis, $G-u\succeq U_{n-1,d}$.

If $v$ lies outside $P(G)$, then $G-u-v\supseteq P_{d+1}$. Thus
$G-u-v\succeq P_{d+1}\succ T_{d,d-2}$.

Suppose that $v$ lies on $P(G)$, then $P(G)$ and $C_r$ have common
vertices, say $v_i,\ldots,v_{i+j}$ with $j\geq 0$.

If $j=0$, i.e., $v=v_i$ is the unique common vertex of $P(G)$ and
$C_r$, then $G-u-v\supseteq P_{i}\cup P_{d-i}\cup P_{2}$. Since
$$m_k(P_{i}\cup P_{d-i}\cup P_{2})\geq m_k(P_{d})\geq m_k(T_{d,d-2})$$
and $m_2(P_{i}\cup P_{d-i}\cup P_{2})>m_2(T_{d,d-2})$, then $P_i\cup
P_{d-i}\cup P_{2}\succ T_{d,d-2}$. Therefore, $G-u-v\succeq
P_{i}\cup P_{d-i}\cup P_{2}\succ T_{d,d-2}$.

If $j>0$. For $v\neq v_i,v_{i+j}$, $G-u-v\supseteq P_{d+1}$. So
$G-u-v\succeq P_{d+1}\succ T_{d,d-2}$. Otherwise, for $v=v_i$ or
$v_{i+j}$, say $v=v_i$. Then $G-u-v\supseteq P_{i}\cup T_1$, where
$T_1\in \mathcal {T}(d-i+1,d-i-1)$ is obtained by attaching a
pendant vertex to vertex $v_{i+j}$ of the path $P=v_{i+1}\cdots
v_d$. For $k\geq 0$,
$$m_k(P_{i}\cup T_1)\geq m_k(P_{i}\cup T_{d-i+1,d-i-1})\geq m_k(T_{d,d-2}).$$
If $(i,j)\neq (1,2)$, then $m_2(P_{i}\cup T_1)>m_2(T_{d,d-2})$,
hence $G-u-v\succeq P_i\cup T_1 \succ T_{d,d-2}$. Otherwise,
$P_i\cup T_1$ is a proper subgraph of $G-u-v$, then $G-u-v\succ
P_i\cup T_1 \succeq T_{d,d-2}$. Thus we always have $G-u-v\succ
T_{d,d-2}$. Therefore, $G-u-v\succ U_{n,d}-u'-v'$.

We have proved that $G-u\succeq U_{n-1,d}$. Then by Lemma
\ref{lem9}, we obtain $G\succ U_{n,d}$.

{\bf{Subcase 2.2}} The neighbor of any pendant vertex outside $P(G)$
also lies outside $C_r$.

If there is a pendant vertex $u$ such that its neighbor $v$ lies
outside $P(G)$, then $G-u-v\supseteq C_r\cup P_{d+1}\supseteq
C_r\cup P_j\cup P_{d-j}$ or $G-u-v\supseteq G'$, where $G' \in
\mathcal {U}(s,d)$ with $d+2\leq s\leq n-2$.

If every pendant vertex outside $P(G)$ is adjacent to a vertex on
$P(G)$, then we choose a pendant vertex $u$, adjacent to $v=v_j$
such that $G-u-v\supseteq C_r\cup P_{j}\cup P_{d-j}$ or
$G-u-v\supseteq P_j\cup G''$, where $G'' \in \mathcal {U}(s',d')$
with $d'\geq d-j-1$, $s'\geq d'+2$ and $s'+j\leq n-2$.

Hence there are three possibilities: $G-u-v\supseteq C_r\cup
P_{j}\cup P_{d-j}$, $G-u-v\supseteq G'$ or $G-u-v\supseteq P_j\cup
G''$.

First, suppose that $G-u-v\supseteq C_r\cup P_{j}\cup P_{d-j}$, then
$$m_k(C_r\cup P_{j}\cup P_{d-j})\geq m_k(P_3\cup P_j\cup P_{d-j})\geq m_k(T_{d,d-2}).$$
In particular, $m_1(C_r\cup P_{j}\cup P_{d-j})> m_1(T_{d,d-2})$.
Thus, $G-u-v \succeq C_r \cup P_j\cup P_{d-j}\succ T_{d,d-2}$.

Next, suppose that $G-u-v\supseteq G'$, then $G-u-v\succeq
U_{s,d}\succeq U_{d+2,d} \succ T_{d+2,d}$.

Finally, suppose that $G-u-v\supseteq P_j\cup G''.$ For $G'' \in
\mathcal {U}(s',d')$ with $s'-d'\leq n-2-j-(d-j-1)=n-d-1$. By the
induction hypothesis, $G''\succeq U_{s',d'}\succeq U_{d-j+1,d-j-1}$.
Thus
$$G-u-v\succeq P_j\cup G''\succeq P_j\cup U_{s',d'}\succeq P_j\cup U_{d-j+1,d-j-1}.$$
For $k\geq 0$,
$$m_k(P_j\cup U_{d-j+1,d-j-1})\geq m_k(P_j\cup T_{d-j+1,d-j-1})\geq m_k(T_{d,d-2}).$$
Furthermore, $m_1(P_j\cup U_{d-j+1,d-j-1})> m_1(T_{d+2,d})$. It
follows that $G-u-v\succeq P_j\cup U_{d-j+1,d-j-1}\succ T_{d,d-2}.$

According to the arguments above, we have proved that $G-u-v\succ
U_{n,d}-u'-v'$. On the other hand, $G-u\succeq U_{n-1,d}$. Thus by
Lemma \ref{lem9}, $G\succ U_{n,d}$.

Combining Cases 1 and 2, we conclude that $G\succ U_{n,d}$ also
holds for $G\in \mathcal {U}(n,d)$ with $3\leq d\leq n-3$ and $G\neq
U_{n,d}$, which yields the result. \qed

\section{Bicyclic graphs with a given diameter}
In what follows we state some new definitions and notations. For a
graph $G\in \mathcal {B}(n)$, it has either two or three distinct
cycles. If $G$ has exactly two cycles, suppose that the lengths of
them are $a$ and $b$ respectively.  If $G$ has   three cycles, then
any two cycles must have at least one edge in common, and we may
choose two cycles of lengths of $a$ and $b$ with $t$ common edges
such that $a-t\geq t$ and $b-t\geq t$. Then, in any case, we choose
two cycles $C_a$ and $C_b$ in $G$. For convenience, let
$C_a=v_0v_1\cdots v_{a-1}v_0$ and $C_b=u_0u_1\cdots u_{b-1}u_0$. If
$C_a$ and $C_b$ have no common edges, then $C_a$ and $C_b$ are
connected by a unique path $P$, say from $v_0$ to $u_0$. Let $l(G)$
be the length of $P$. If $C_a$ and $C_b$ have exactly $t(\geq1)$
common edges, and thus have exactly $t+1$ common vertices, say,
$v_0=u_0, v_1=u_1, \ldots, v_t=u_t$, then $C_c=u_0u_{b-1}\cdots
u_{t+1}u_tv_{t+1}v_{t+2}\cdots v_{a-1}v_0$ is the third cycle of
$G$, where $c=b+a-2t$. If we write $w_0=u_0, w_1=u_{b-1},\ldots,
w_{c-1}=v_{a-1}$, then $C_c=w_0w_1\cdots w_{c-1}w_0$. Denote by
$d(G)$ the diameter of $G$.

Now we turn our attention to the minimal matching energy of graphs
in $\mathcal {B}(n,d)$ with $3\leq d\leq n-2$. We first deal with
the case $d=n-3$.

\begin{lemma}\label{lem $n-d=3$}
Let $G\in \mathcal {B}(n,n-3)$ with $n\geq 7$, and $G \neq
B_{n,n-3}$. Then $G\succ B_{n,n-3}$.
\end{lemma}
\pf By induction on $n$  to prove this fact.

For $n=7$ and $n=8$, there are only finitely many graphs we need to
consider. Then by Lemma \ref{lem1} and direct check, we can get
$G\succ B_{n,n-3}$.

Suppose that the result holds for all graphs in $\mathcal
{B}(n-1,n-4)$ and $\mathcal {B}(n-2,n-5)$, where $n\geq 9$. Let
$G\in \mathcal {B}(n,n-3)$ and $G \neq B_{n,n-3}$.

{\bf{Case 1}} There is a pendent vertex $u$ in $G$ such that the
degree of its neighbor $v$ is $2$.

In this case, $G-u\in \mathcal {B}(n-1,n-4)$ and $G-u-v\in \mathcal
{B}(n-2,n-5)$. Since $G \neq B_{n,n-3}$, then $G-u \neq B_{n-1,n-4}$
or $G-u-v \neq B_{n-2,n-5}$. By the induction hypothesis, $G-u \succ
B_{n-1,n-4}$ and $G-u-v \succeq B_{n-2,n-5}$, or $G-u \succeq
B_{n-1,n-4}$ and $G-u-v \succ B_{n-2,n-5}$. Hence, $G\succ
B_{n,n-3}$.

{\bf{Case 2}} The neighbor of any pendent vertex has degree at least
$3$ or there is no pendent vertex.

\begin{figure}[h,t,b,p]
\begin{center}
\includegraphics[scale = 0.9]{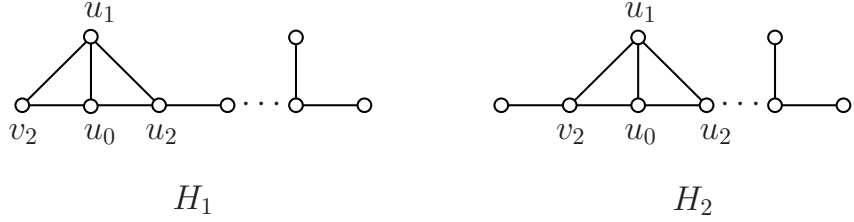}
\caption{The graphs $H_i$ for $i=1,2$.}\label{H_i}
\end{center}
\end{figure}

Then $G$ is isomorphic to some $H_j$, $j=1, 2$(see Figure
\ref{H_i}), or $G$ contains one triangle or one quadrangle which has
at most one common vertex with the other cycle that is a triangle or
a quadrangle.

\begin{figure}[h,t,b,p]
\begin{center}
\includegraphics[scale = 0.8]{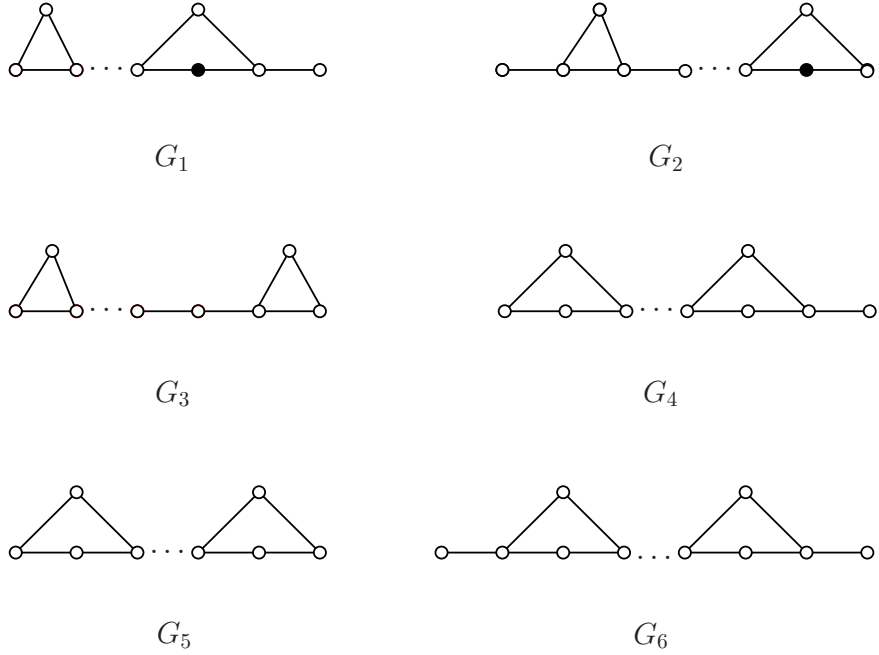}
\caption{The graphs $G_i$ for $i=1,2,3,4,5,6$.}\label{G_i}
\end{center}
\end{figure}

If $G$ is isomorphic to $H_1$, then by Lemmas \ref{lem1}, \ref{TT},
\ref{TdTd_0} and Theorem \ref{th2},
\begin{eqnarray*}
m_k(G) &=& m_k(G-u_1v_2)+m_{k-1}(G-u_1-v_2)\\
       &\geq& m_k(U_{n,n-3})+m_{k-1}(T_{n-2,d-1})\\
       &\geq& m_k(U_{n,n-3})+m_{k-1}(P_{n-6}\cup S_4)\\
       &=&m_k(B_{n,n-3}).
\end{eqnarray*}
Moreover, $m_2(G)>m_2(B_{n,n-3})$, thus $G\succ B_{n,n-3}$.

If $G$ is isomorphic to $H_2$, then by Lemmas \ref{lem1}, \ref{PPP}
and Theorem \ref{th2},
\begin{eqnarray*}
m_k(G) &=& m_k(G-u_1u_0)+m_{k-1}(G-u_1-u_0)\\
       &\geq& m_k(U_{n,n-3})+m_{k-1}(P_2\cup T_{n-4,n-6})\\
       &\geq& m_k(U_{n,n-3})+m_{k-1}(P_2\cup P_{n-5})+m_{k-2}(P_2\cup P_{n-7})\\
       &\geq& m_k(U_{n,n-3})+m_{k-1}(P_3\cup P_{n-6})+m_{k-2}(P_{n-6})\\
       &=&m_k(B_{n,n-3}).
\end{eqnarray*}
Similarly, $m_2(G)>m_2(B_{n,n-3})$, thus $G\succ B_{n,n-3}$.

Otherwise, $G$ contains one triangle or one quadrangle which has at
most one common vertex with the other cycle that is a triangle or a
quadrangle. Choose $C_a$ and $C_b$ as above. Let $b\geq a$.

If $a=3$, then $G$ is isomorphic to $G_1$, $G_2$ or $G_3$ in Figure
\ref{G_i}, where the black vertices may not occur. Similarly, we can
obtain that $G\succ B_{n,n-3}$.

If $a=4$, then $G$ is isomorphic to $G_4$, $G_5$ or $G_6$ in Figure
\ref{G_i}. We can show that $G\succ B_{n,n-3}$ in the same way.

Hence the conclusion follows. \qed

\begin{lemma}\label{GB_{n,d+1}}
Let $G\in \mathcal {B}(n,d)$ with $n\geq 8$ and $3\leq d\leq n-4$.
If $G$ contains no pendent vertices, then $G\succ B_{n,d+1}$.
\end{lemma}
\pf We choose $C_a$, $C_b$ in $G$ and if there exists the third
cycle, then we choose $C_c$ and $t$ as above. Let $b\geq a$. Since
$d\leq n-4$, we have $b\geq5$.

{\bf{Case 1}} $C_a$ and $C_b$ have no common edges.

Then $d=\lfloor \frac{a}{2} \rfloor+\lfloor \frac{b}{2}
\rfloor+l(G)$, $d(G-u_1u_2)=\lfloor \frac{a}{2} \rfloor+l(G)+b-2\geq
d+1$, $d(G-u_1-u_2)=\lfloor \frac{a}{2} \rfloor+l(G)+b-3\geq d$.
According to Lemmas \ref{lem1}, \ref{lem2}, \ref{Udd_0} and Theorem
\ref{th2},
\begin{eqnarray*}
m_k(G) &=& m_k(G-u_1u_2)+m_{k-1}(G-u_1-u_2)\\
       &\geq& m_k(U_{n,d+1})+m_{k-1}(U_{n-2,d})\\
       &\geq& m_k(U_{n,d+1})+m_{k-1}(P_{d-2}\cup S_{n-d})\\
       &=&m_k(B_{n,d+1}).
\end{eqnarray*}

Further, we have $m_2(G)>m_2(B_{n,d+1})$, thus $G\succ B_{n,d+1}$.

{\bf{Case 2}} $C_a$ and $C_b$ have at least one common edge.

Notice that $a-t\geq t, b-t\geq t$, where $t\geq 1$. It follows that
$c\geq b+1$, $d=\lfloor \frac{c}{2}\rfloor=\lfloor
\frac{(a+b)}{2}\rfloor-t$, $d(G-w_0-w_1)=c-3\geq d$ and $d\geq 3$.

If $b>5$ or $b=5$ and $a$ is even, then $d(G-w_0w_1)=\lfloor
\frac{a}{2} \rfloor+b-t-1\geq d+1$, by Lemmas \ref{lem1},
\ref{lem2}, \ref{BUT} and Theorem \ref{th2},
\begin{eqnarray*}
m_k(G) &=& m_k(G-w_0w_1)+m_{k-1}(G-w_0-w_1)\\
       &\geq& m_k(U_{n,d+1})+m_{k-1}(U_{n-2,d})\\
       &\geq& m_k(U_{n,d+1})+m_{k-1}(T_{n-2,d})\\
       &\geq& m_k(U_{n,d+1})+m_{k-1}(P_{d-2}\cup S_{n-d})\\
       &=&m_k(B_{n,d+1})
\end{eqnarray*}
together with $m_2(G)>m_2(B_{n,d+1})$, hence $G\succ B_{n,d+1}$.

\begin{figure}[h,t,b,p]
\begin{center}
\includegraphics[scale = 0.8]{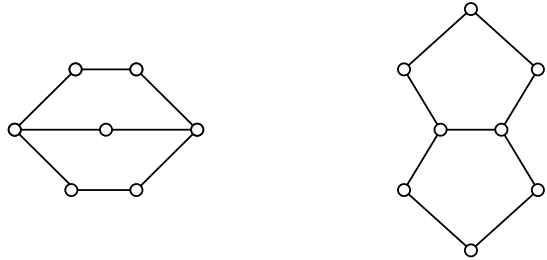}
\caption{The graphs isomorphic to $G$ when $b=5$ and $a$ is
odd.}\label{b=5}
\end{center}
\end{figure}

If $b=5$ and $a$ is odd, then $G$ is isomorphic to the two graphs in
Figure \ref{b=5}, it is easy to verify  that $G\succ B_{n,d+1}$ and
the proof is complete. \qed

\begin{lemma}\label{B_{n,d+1}}
Let $G\in \mathcal {B}(n,d)$ with $n\geq 8$ and $3\leq d\leq n-4$.
If $G$ contains exactly one pendent vertex $u$ on all diametrical
paths of $G$ such that $G-u$ contains no pendent vertices, then
$G\succ B_{n,d+1}$.
\end{lemma}
\pf We choose $C_a$, $C_b$ in $G$ and if there exists the third
cycle, then we choose $C_c$ and $t$ as above. Let $b\geq a$. Since
$d\leq n-4$, we have $b\geq5$. Let $v$ be the neighbor of $u$.

{\bf{Case 1}} $C_a$ and $C_b$ have no common edges.

Then $d=\lfloor \frac{a}{2} \rfloor+\lfloor \frac{b}{2}
\rfloor+l(G)+1$.

If $b\geq 7$, then $d(G-u_1u_2)\geq \lfloor \frac{a}{2}
\rfloor+l(G)+b-2\geq d+1$.

If $v$ lies on $C_a$ and $b=5$ or $6$, then $d(G-u_1u_2)=\lfloor
\frac{a}{2}\rfloor+l(G)+b-1=d+1$.

If $v$ lies on $C_b$, $a=5$ and $b=6$, then $d(G-v_1v_2)=\lfloor
\frac{b}{2}\rfloor+l(G)+a-1=d+1$.

In these cases, the proof is the same as Case 1 of Lemma
\ref{GB_{n,d+1}}.

Otherwise, $v$ lies on $C_b$, $a=3$ or $4$ and $b=5$ or $6$.

If $l(G)=0$, then $G$ is isomorphic to finitely many graphs. Apply
Lemma \ref{lem1} and direct calculation, we can get $G\succ
B_{n,d+1}$.

So suppose that $l(G)\geq1$. If $a=3$ and $b=5$, then $G-v_0v_1\in
\mathcal {U}(n,d+1)$. $G-v_0-v_1=P_1\cup G_0$, where $G_0\in
\mathcal {U}(n-3,d-2)$. Since $d(G_0)=d-2$, $n(G_0)=n-3$, then
$n-3=d-2+3$, that is, $n-d=4$, meanwhile, $l(G)-1=d-2-3=d-5$, i.e.,
$l(G)=d-4$. Thus $G_0-u_0u_1\supseteq P_{d-2}\cup S_3,
G_0-u_0-u_1=P_{d-5}\cup P_4$. Therefore, by Lemmas \ref{lem1},
\ref{lem2}, \ref{PPP} and Theorem \ref{th2},
\begin{eqnarray*}
m_k(G) &=& m_k(G-v_0v_1)+m_{k-1}(G-v_0-v_1)\\
       &\geq& m_k(U_{n,d+1})+m_{k-1}(G_0-u_0u_1)+m_{k-2}(G_0-u_0-u_1)\\
       &\geq& m_k(U_{n,d+1})+m_{k-1}(P_{d-2}\cup S_3)+m_{k-2}(P_{d-2})\\
       &\geq& m_k(U_{n,d+1})+m_{k-1}(P_{d-2}\cup S_{n-d})\\
       &=&m_k(B_{n,d+1}).
\end{eqnarray*}
In particular, $m_2(G)>m_2(B_{n,d+1})$, thus $G\succ B_{n,d+1}$.

If $a=3$ and $b=6$, by similar arguments, we can obtain that $G\succ
B_{n,d+1}$.

If $a=4$ and $b=5$, then $n=d+1+3=d+4$, i.e., $n-d-1=3$,
$G-v_0v_1\in \mathcal {U}(n,d+1)$ and $G-v_0-v_1=P_2\cup G_0$, where
$G_0\in \mathcal {U}(n-4,d-3)$. Similarly,
\begin{eqnarray*}
m_k(G) &=& m_k(G-v_0v_1)+m_{k-1}(G-v_0-v_1)\\
       &\geq& m_k(U_{n,d+1})+m_{k-1}(P_2\cup G_0-u_0u_1)+m_{k-2}(P_2\cup G_0-u_0-u_1)\\
       &\geq& m_k(U_{n,d+1})+m_{k-1}(P_2\cup P_{d-2})+m_{k-2}(P_2\cup P_{d-6}\cup P_4)\\
       &\geq& m_k(U_{n,d+1})+m_{k-1}(P_3\cup P_{d-3})+m_{k-2}(P_{d-2})\\
       &=&m_k(B_{n,d+1}).
\end{eqnarray*}
Moreover, $m_2(G)>m_2(B_{n,d+1})$, hence $G\succ B_{n,d+1}$.

If $a=4$ and $b=6$, we can verify that $G\succ B_{n,d+1}$ in the
same way.

{\bf{Case 2}} $C_a$ and $C_b$ have at least one common edge.

Then $d=\lfloor \frac{c}{2} \rfloor+1=\lfloor \frac{(a+b)}{2}
\rfloor-t+1$. Since $b\geq 5$, $w_0$, $w_1\neq v$ and $d(G-w_0w_1)=
\lfloor \frac{a}{2} \rfloor+b-t-1$. When $b\geq 6$, $d(G-w_0w_1)\geq
d+1$, $d(G-w_0w_1)\geq d-1$. And then
\begin{eqnarray*}
m_k(G) &=& m_k(G-w_0w_1)+m_{k-1}(G-w_0-w_1)\\
       &\geq& m_k(U_{n,d+1})+m_{k-1}(T_{n-2,d-1})\\
       &\geq& m_k(U_{n,d+1})+m_{k-1}(P_{d-2}\cup S_{n-d})\\
       &=&m_k(B_{n,d+1})
\end{eqnarray*}
along with $m_2(G)>m_2(B_{n,d+1})$, hence $G\succ B_{n,d+1}$. Now,
we are left with the cases:

(i) $b=5, t=2, a=4, uv\in C_b$;

(ii) $b=5, t=2, a=4, uv\in C_a$;

(iii) $b=5, t=2, a=5$.

It can be checked directly that $G\succ B_{n,d+1}$ in these cases.

Combining Cases 1 and 2, we arrive at the result. \qed

\begin{theorem}\label{thm3}
Let $G\in \mathcal {B}(n,d)$ with $n\geq 8$ and $3\leq d\leq n-3$.
If there are two vertex-disjoint cycles in $G$, then $G\succ
B_{n,d}$.
\end{theorem}
\pf The proof proceeds by induction on $n-d$. By Lemma \ref{lem
$n-d=3$}, the result holds for $n-d=3$. Let $h\geq 4$ and assume
that the result holds for $n-d<h$. Suppose that $G\in \mathcal
{B}(n,d)$ with $n-d=h$.

{\bf{Case 1}} There is no pendent vertex in $G$.

Then by Lemmas \ref{Bdd_0} and \ref{GB_{n,d+1}}, $G\succ
B_{n,d+1}\succ B_{n,d}$.

{\bf{Case 2}} There is a pendent vertex outside some diametrical
path $P(G)=x_0x_1\cdots x_d$.

Let $u$, adjacent to $v$, be a pendent vertex outside $P(G)$ in $G$.
Then $G-u\in \mathcal {B}(n-1,d)$. Since $(n-1)-d<h$, by the
induction hypothesis, $G-u\succ B_{n-1,d}$.

By Lemma \ref{lem1},
$$m_k(B_{n,d})=m_k(B_{n-1,d})+m_{k-1}(T_{d+1,d-2}) \eqno(*)$$
Meanwhile, let $H=G-u-v$, then
\begin{eqnarray*}
m_k(G) &=& m_k(G-u)+m_{k-1}(G-u-v)\\
       &=& m_k(G-u)+m_{k-1}(H)\\
       &\geq& m_k(B_{n-1,d})+m_{k-1}(H).
\end{eqnarray*}
Hence, to complete the proof we shall show that $m_k(H)\geq
m_k(T_{d+1,d-2})$. Select $C_a$ and $C_b$ as above.

{\bf{Subcase 2.1}} $v$ lies on some cycle, say $C_a$. Then
$H\supseteq C_b$.

First, suppose that $P(G)$ and $C_b$ have no common vertices. Then
$H\supseteq P_i \cup P_{d-i}\cup C_b$ when $v$ lies on $P(G)$, say
$v=v_i$. Otherwise, $H\supseteq P_{d+1}\cup C_b\supseteq P_i\cup
P_{d-i}\cup C_b$. Thus
\begin{eqnarray*}
m_k(H) &\geq& m_k(P_i \cup P_{d-i}\cup C_b)\geq m_k(P_{d-1}\cup S_b)\geq m_k(P_{d-1}\cup P_3)\\
       &\geq& m_k(P_{d+1})\geq m_k(T_{d+1,d-2}).
\end{eqnarray*}

Next, suppose that $P(G)$ and $C_b$ have common vertices
$x_l,\ldots,x_{l+q}$, where $q\geq0$.

If $v$ lies outside $P(G)$, then $H\supseteq G_1$, where $G_1\in
\mathcal {U}(s_1,d)$ with $s_1\geq d+2$. Hence,
$$m_k(H)\geq m_k(G_1) \geq m_k(U_{s_1,d})\geq m_k(T_{s_1,d})\geq m_k(P_{d+1})\geq m_k(T_{d+1,d-2}).$$

So suppose that $v$ lies on $P(G)$. Then $P(G)$ and $C_a$ have
common vertices, say $x_i,\ldots,x_{i+p}$, where $p\geq0$.

When $p=0$, then $i\geq1$, thus $H\supseteq P_2 \cup P_{i}\cup G_2$,
where $G_2\in \mathcal {U}(s_2,d_2)$, $d_2\geq d-i-1\geq1$ and
$s_2\geq d_2+2$. If $d_2=1, i=d-2$, then $G_2=C_3$. Therefore,
$$m_k(H)\geq m_k(P_2 \cup P_{d-2}\cup C_3) \geq m_k(P_{d-1}\cup P_3)\geq m_k(P_{d+1})\geq m_k(T_{d+1,d-2}).$$
If $d_2=2$, then $s_2\geq 4, i\geq d-3$. Consequently,
\begin{eqnarray*}
m_k(H) &\geq&m_k(P_2 \cup P_{i}\cup G_2)\geq m_k(P_2 \cup P_{i}\cup S_{s_2})\geq m_k(P_2 \cup P_{i}\cup T_{4,2})\\
       &\geq& m_k(T_{i+4,i+2})\geq m_k(T_{d+1,d-1})\geq m_k(T_{d+1,d-2}).
\end{eqnarray*}
If $d_2\geq3$, then $s_2\geq 5$.  Thus
\begin{eqnarray*}
m_k(H) &\geq& m_k(P_2 \cup P_{i}\cup G_2)\geq m_k(P_2 \cup P_{i}\cup U_{s_2,d_2})\geq m_k(P_2 \cup P_{i}\cup T_{s_2,d_2})\\
       &\geq& m_k(T_{s_2+i,d_2+i})\geq m_k(T_{d_2+i+2,d_2+i})\geq m_k(T_{d+1,d-1})\geq m_k(T_{d+1,d-2}).
\end{eqnarray*}

When $p=1$. If $v=v_i$, then $i\geq1$ and $H\supseteq P_i\cup G_3$,
where $G_3\in \mathcal {U}(s_3,d_3), d_3\geq d-i \geq 3, s_3\geq
d_3+2$. Accordingly,
\begin{eqnarray*}
m_k(H) &\geq& m_k(P_{i}\cup G_3)\geq m_k(P_{i}\cup U_{s_3,d_3})\geq m_k(P_{i}\cup T_{s_3,d_3})\geq m_k(T_{s_3+i-1,d_3+i-1})\\
&\geq& m_k(T_{d_2+i+2,d_2+i})\geq m_k(T_{d+1,d-1})\geq
m_k(T_{d+1,d-2}).
\end{eqnarray*}
If $v=v_{i+1}$, then $H\supseteq P_{i+2}\cup G_4$, where $G_4\in \mathcal {U}(s_4,d_4), d_4\geq d-(i+1)-1=d-i-2\geq1, s_4\geq d_4+2$.\\
If $d_4=1$, then $i=d-3, G_4=C_3$, we have
$$m_k(H)\geq m_k(P_{d-1}\cup C_3) \geq m_k(P_{d-1}\cup P_3)\geq m_k(P_{d+1})\geq m_k(T_{d+1,d-2}).$$
If $d_4=2$, then $s_4\geq 4, i\geq d-4$. Consequently,
\begin{eqnarray*}
m_k(H) &\geq& m_k(P_{i+2}\cup G_4)\geq m_k(P_{i+2}\cup S_{s_4})\geq m_k(P_{i+2}\cup T_{4,2})\\
       &\geq& m_k(T_{i+5,i+3})\geq m_k(T_{d+1,d-1})\geq m_k(T_{d+1,d-2}).
\end{eqnarray*}
If $d_4\geq3$. Then
\begin{eqnarray*}
m_k(H) &\geq& m_k(P_{i+2}\cup G_4)\geq m_k(P_{i+2}\cup U_{s_4,d_4})\geq m_k(P_{i+2}\cup T_{s_4,d_4})\\
       &\geq& m_k(T_{s_4+i+1,d_4+i+1})\geq m_k(T_{d+1,d-1})\geq m_k(T_{d+1,d-2}).
\end{eqnarray*}

Now suppose that $p\geq2$. If $v\neq x_i$, $x_{i+p}$, then
$H\supseteq G_5$, where $G_5\in \mathcal {U}(s_5,d_5), d_5\geq d,
s_5\geq d_5+2$. Hence
$$m_k(H)\geq m_k(G_5) \geq m_k(U_{s_5,d_5})\geq m_k(T_{s_5,d_5})\geq m_k(T_{d+1,d-2}).$$
If $v=x_i$, then $H\supseteq P_i\cup G_6$, where $G_6\in \mathcal
{U}(s_6,d_6), d_6\geq d-i-1\geq3, s_6\geq d_6+3$. Therefore,
\begin{eqnarray*}
m_k(H)&\geq& m_k(P_i\cup G_6) \geq m_k(P_i\cup U_{s_6,d_6})\geq m_k(P_i\cup T_{s_6,d_6})\\
&\geq& m_k(T_{s_6+i-1,d_6+i-1}) \geq m_k(T_{d+1,d-2}).
\end{eqnarray*}
If $v=x_{i+p}$, then $H\supseteq T_1\cup G_7$ or $P_{i+p+1}\cup
G_7$, where $G_7\in \mathcal {U}(s_7,d_7), d_7\geq d-i-p-1\geq1,
s_7\geq d_7+2, T_1 \in \mathcal {T}(i+p+1,d_{T_{1}}), d_{T_{1}}\geq
i+p-1$. If $d_7=1$, then $i+p=d-2, G_7=C_3$, thus
$$m_k(H)\geq m_k(T_{d-1,d-3}\cup C_3)\geq m_k(T_{d-1,d-3}\cup P_3)\geq m_k(T_{d+1,d-2})$$
or
\begin{eqnarray*}
m_k(H) &\geq& m_k(P_{i+p+1}\cup G_7)\geq m_k(P_{i+p+1}\cup C_3)\geq m_k(P_{i+p+1}\cup P_3)\\
       &\geq& m_k(P_{i+p+3})=m_k(P_{d+1})\geq m_k(T_{d+1,d-2}).
\end{eqnarray*}
If $d_7=2$, then $i+p\geq d-3, s_7\geq 4$, accordingly,
\begin{eqnarray*}
m_k(H) &\geq& m_k(T_{i+p+1,i+p-1}\cup S_{s_7})\geq m_k(T_{i+p+1,i+p-1}\cup T_{4,2})\\
       &\geq& m_k(T_{i+p+4,i+p+1})\geq m_k(T_{d+1,d-2})
\end{eqnarray*}
or
$$m_k(H)\geq m_k(P_{i+p+1}\cup S_{s_7})\geq m_k(P_{i+p+1}\cup T_{4,2})\geq m_k(T_{i+p+4,i+p+2}) \geq m_k(T_{d+1,d-2}).$$
If $d_7\geq3$, then
\begin{eqnarray*}
m_k(H) &\geq& m_k(T_{i+p+1,i+p-1}\cup U_{s_7,d_7})\geq m_k(T_{i+p+1,i+p-1}\cup T_{s_7,d_7})\\
       &\geq& m_k(T_{i+p+s_7,i+p+d_7-1})\geq m_k(T_{d+1,d-2})
\end{eqnarray*}
or
$$m_k(H)\geq m_k(P_{i+p+1}\cup U_{s_7,d_7})\geq m_k(P_{i+p+1}\cup T_{s_7,d_7})\geq m_k(T_{i+p+s_7,i+p+d_7})\geq m_k(T_{d+1,d-2}).$$

{\bf{Subcase 2.2}} $v$ lies outside any cycle. Then $H\supseteq
C_a\cup C_b$.

First, suppose that $v$ lies on $P(G)$ and take $v=x_i$. If $P(G)$
has no common vertices with any cycle, then $H\supseteq C_a\cup
C_b\cup P_i\cup P_{d-i}$. Thus
$$m_k(H)\geq m_k(C_a\cup C_b\cup P_i\cup P_{d-i})\geq m_k(P_3\cup P_{d-1})\geq m_k(P_{d+1})\geq m_k(T_{d+1,d-2}).$$
If $P(G)$ has no common vertices with exactly one cycle, say $C_a$.
Then $H\supseteq C_a\cup P_i\cup G_1$, where $G_1\in \mathcal
{U}(s_1,d_1), d_1\geq d-i-1, s_1\geq d_1+2$. If $d_1=1$, then
$i=d-2, G_1=C_3$, hence
$$m_k(H)\geq m_k(C_a\cup P_i\cup C_3)\geq m_k(P_2\cup P_i\cup P_3)\geq m_k(P_{i+3})=m_k(P_{d+1})\geq m_k(T_{d+1,d-2}).$$
If $d_1=2$, then $i\geq d-3, s_1\geq 4$. Consequently,
$$m_k(H)\geq m_k(C_a\cup P_i\cup G_1)\geq m_k(P_2\cup P_i \cup T_{s_1,2})\geq m_k(T_{i+s_1,i+2})\geq m_k(T_{d+1,d-2}).$$
If $d_1\geq 3$, then
\begin{eqnarray*}
m_k(H) &\geq& m_k(C_a\cup P_i\cup G_1)\geq m_k(P_2\cup P_i \cup U_{s_1,d_1})\geq m_k(P_2\cup P_i \cup T_{s_1,d_1})\\
       &\geq& m_k(P_{i+1}\cup T_{s_1,d_1})\geq m_k(T_{s_1+i,d_1+i})\geq m_k(T_{d+1,d-2}).
\end{eqnarray*}
If $P(G)$ has common vertices with both cycles, then $H\supseteq P_i\cup G_2$ or $G_3\cup G_4$, where $G_2\in \mathcal {U}(s_2,d_2)$, $G_3\in \mathcal {U}(s_3,d_3)$, $G_4\in \mathcal {U}(s_4,d_4)$. Meanwhile, $d_2\geq d-i-1\geq3$, $n-2-i\geq s_2\geq d_2+3$, $d_3\geq i-1\geq1$, $s_3\geq d_3+2$, $d_4\geq d-i-1\geq1$, $s_4\geq d_4+2$.\\
Suppose that $H\supseteq P_i\cup G_2$. Since $s_2-d_2\leq
n-2-i-(d-i-1)=n-d-1<h$, thus by the induction hypothesis, $G_2\succ
B_{s_2,d_2}$. Then
\begin{eqnarray*}
m_k(H)&\geq& m_k(P_i\cup G_2)\geq m_k(P_i\cup B_{s_2,d_2})\geq m_k(P_i\cup T_{s_2,d_2})\\
&\geq& m_k(T_{s_2+i-1,d_2+i-1})\geq m_k(T_{d+1,d-2}).
\end{eqnarray*}
Suppose that $H\supseteq G_3\cup G_4$. If $d_3=d_4=1$, then $n=8$, $d=4$, $G_3=G_4=C_3$. In this case, it is easy to obtain that $G\succ B_{8,4}$.\\
If $d_3=2$, $d_4=1$, then $d=5$, $s_3\geq 4$, $G_4=C_3$. We can have
$$m_k(H)\geq m_k(G_3\cup G_4)\geq m_k(S_{s_3}\cup C_3)\geq m_k(T_{4,2}\cup P_3)\geq m_k(T_{6,4})\geq m_k(T_{d+1,d-2}).$$
If $d_3\geq 3$, $d_4=1$, then $d_3\geq d-3$. Accordingly,
\begin{eqnarray*}
m_k(H)&\geq& m_k(G_3\cup G_4)\geq m_k(U_{s_3,d_3}\cup C_3)\geq m_k(T_{s_3,d_3}\cup P_3)\\
&\geq& m_k(T_{s_3+2,d_3+2})\geq m_k(T_{d+1,d-2}).
\end{eqnarray*}
If $d_3=2$, $d_4=2$, then $d=6$, $s_3\geq 4$, $s_4\geq 4$. Hence
$$m_k(H)\geq m_k(G_3\cup G_4)\geq m_k(S_{s_3}\cup S_{s_4})\geq m_k(T_{4,2}\cup T_{4,2})\geq m_k(T_{7,4})=m_k(T_{d+1,d-2}).$$
If $d_3\geq 3$, $d_4=2$, then $d_3\geq d-4$, $s_4\geq 4$. Thus
\begin{eqnarray*}
m_k(H)&\geq& m_k(G_3\cup G_4)\geq m_k(U_{s_3,d_3}\cup S_{s_4})\geq m_k(T_{s_3,d_3}\cup T_{s_4,2})\\
&\geq& m_k(T_{s_3+s_4-1,d_3+2})\geq m_k(T_{d+1,d-2}).
\end{eqnarray*}
If $d_3\geq 3$, $d_4\geq3$, then $d_3+d_4\geq d-2$. Therefore,
\begin{eqnarray*}
m_k(H)&\geq& m_k(G_3\cup G_4)\geq m_k(U_{s_3,d_3}\cup U_{s_4,d_4})\geq m_k(T_{s_3,d_3}\cup T_{s_4,d_4})\\
&\geq& m_k(T_{s_3+s_4-1,d_3+d_4})\geq m_k(T_{d+1,d-2}).
\end{eqnarray*}

Next, suppose that $v$ lies outside $P(G)$. Then $H\supseteq C_a\cup
C_b\cup P(G)$, $C_a\cup G_5$ or $G_6$, where $G_5\in \mathcal
{U}(s_5,d)$ with $s_5\geq d+2$ and $G_6\in \mathcal {B}(s_6,d)$ with
$d+3\leq s_6\leq n-2$. It is easy to show as above that $m_k(H)\geq
m_k(T_{d+1,d-2})$.

{\bf{Case 3}} Any diametrical path of $G$ contains all pendent
vertices.

Let $P(G)=x_0x_1\cdots x_d$ be any diametrical path of $G$. Suppose
that $y_0y_1\cdots y_p$ is a path whose internal vertices
$y_1,y_2,\ldots,y_{p-1}$ all have degree two and $y_p$ is a pendent
vertex. Then we call it a pendent path, denoted by $(y_0,y_p)$.

{\bf{Subcase 3.1}} There are exactly two pendent vertices in $G$,
namely, $x_0$ and $x_d$.

Suppose that $deg_G(x_i)$, $deg_G(x_l)\geq 3$ such that $(x_i,x_0)$ and $(x_l,x_d)$ are distinct pendent paths. Let $s=l-i$.\\
If $s=0$, i.e., $x_i=x_l$. Then $i\geq3, l\leq d-3$. Since
\begin{eqnarray*}
m_k(G) &=& m_k(G-x_{i-3}x_{i-2})+m_{k-1}(G-x_{i-3}-x_{i-2})\\
       &=& m_k(G-x_{i-3}x_{i-2}-x_{l+1}x_{l+2})+m_{k-1}(G-x_{i-3}x_{i-2}-x_{l+1}-x_{l+2})\\
       && +~ m_{k-1}(G-x_{i-3}-x_{i-2}-x_{l+2}x_{l+3})+m_{k-2}(G-x_{i-3}-x_{i-2}-x_{l+2}-x_{l+3})\\
       &=& m_k(G_1\cup P_{i-2}\cup P_{d-i-1})+m_{k-1}(G_3\cup P_{i-2}\cup P_{d-i-2})\\
       && +~ m_{k-1}(G_2\cup P_{i-3}\cup P_{d-i-2})+m_{k-2}(G_4\cup P_{i-3}\cup P_{d-i-3})
\end{eqnarray*}
and
\begin{eqnarray*}
m_k(B_{n,d}) &=& m_k(B_{n-d+3,3}\cup P_{d-3})+m_{k-1}(S_{n-d+2}\cup P_{d-4})\\
       &=& m_k(B_{n-d+3,3}\cup P_{i-2}\cup P_{d-i-1})+m_{k-1}(B_{n-d+3,3}\cup P_{i-3}\cup P_{d-i-2})\\
&& +~ m_{k-1}(S_{n-d+2}\cup P_{i-2}\cup
P_{d-i-2})+m_{k-2}(S_{n-d+2}\cup P_{i-3}\cup P_{d-i-3}),
\end{eqnarray*}
it suffices to prove that $G_1,G_2 \succ B_{n-d+3,3}$ and $G_3,
G_4\succ S_{n-d+2}$, where $G_1=G-(x_{i-3},x_0)-(x_{l+2},x_d)\in
\mathcal {B}(n-d+3,d_1)$, $G_2=G-(x_{i-2},x_0)-(x_{l+3},x_d)\in
\mathcal {B}(n-d+3,d_2)$, $G_3=G-(x_{i-3},x_0)-(x_{l+1},x_d)$,
$G_4=G-(x_{i-2},x_0)-(x_{l+2},x_d)$, $d_1\geq 4$, $d_2\geq 4$. Since
$n-d+3-d_1\leq n-d-1<h$, $n-d+3-d_2\leq n-d-1<h$. Then by the
induction hypothesis, $G_1\succ B_{n-d+3,d_1}\succ B_{n-d+3,3}$,
$G_2\succ B_{n-d+3,d_2}\succ B_{n-d+3,3}$.
In addition, both $G_3$ and $G_4$ are bicyclic graphs with $n-d+2$ vertices, consequently, we have $G_3,G_4\succ S_{n-d+2}$.\\
If $s=1$ or $s=2$, then by similar arguments as above, we have the desired result.\\
If $s\geq 3$, it is easy to obtain that $i\geq 2$ and $l\leq d-2$.
Then
\begin{eqnarray*}
m_k(G) &=& m_k(G-x_{i-2}x_{i-1})+m_{k-1}(G-x_{i-2}-x_{i-1})\\
       &=& m_k(G-x_{i-2}x_{i-1}-x_{l+1}x_{l+2})+m_{k-1}(G-x_{i-2}x_{i-1}-x_{l+1}-x_{l+2})\\
       &&+~ m_{k-1}(G-x_{i-2}-x_{i-1}-x_{l+1}x_{l+2})+m_{k-2}(G-x_{i-2}-x_{i-1}-x_{l+1}-x_{l+2})\\
       &=& m_k(G_5\cup P_{i-1}\cup P_{d-l-1})+m_{k-1}(G_7\cup P_{i-1}\cup P_{d-l-2})\\
       &&+~ m_{k-1}(G_6\cup P_{i-2}\cup P_{d-l-1})+m_{k-2}(G_8\cup P_{i-2}\cup P_{d-l-2}),
\end{eqnarray*}
\begin{eqnarray*}
m_k(B_{n,d}) &=& m_k(B_{n-d+l+1,l+1}\cup P_{d-l-1})+m_{k-1}(B_{n-d+l,l}\cup P_{d-l-2})\\
       &=& m_k(B_{n-d+s+2,s+2}\cup P_{i-1}\cup P_{d-l-1})+m_{k-1}(B_{n-d+s+1,s+1}\cup P_{i-2}\cup P_{d-l-1})\\
       &&+~ m_{k-1}(B_{n-d+s+1,s+1}\cup P_{i-1}\cup P_{d-l-2})+m_{k-2}(B_{n-d+s,s}\cup P_{i-2}\cup P_{d-l-2}).
\end{eqnarray*}
Hence it suffices to show that $G_5\succ B_{n-d+s+2,s+2}$, $G_6,
G_7\succ B_{n-d+s+1,s+1}$, $G_8\succ B_{n-d+s,s}$.

Let $d_j=d(G_j)$, and $n_j=|V(G_j)|$, where $j=5, 6, 7, 8$. Then
$d_j\geq4$. If $n_j-d_j<h$ holds for all $j\in \{5, 6, 7, 8\}$, then
by the induction hypothesis and previous Lemmas, we have the desired
results. Otherwise, there exists at least a $j\in \{5, 6, 7, 8\}$
such that $n_j-d_j=h$. When $j=5$, $G_5\in \mathcal
{B}(n-d+s+2,s+2)$. If there exists some diametrical path $P(G_5)$
such that $x_{i-1}$ or $x_{l+1}$ lies outside $P(G_5)$, the proof is
similar with Case 2, thus $G_5\succ B_{n-d+s+2,s+2}$. Otherwise,
$G_5-x_{i-1}\in \mathcal {B}(n-d+s+1,s+1)$, then by Lemma
\ref{B_{n,d+1}}, $G_5\succ B_{n-d+s+1,s+2}$. We also have
$G_5-x_{i-1}-x_i\succeq U_{n-d+s,s}\succ T_{n-d+s,s}\succ
T_{s+3,s}$. Therefore, $G_5\succ B_{n-d+s+2,s+2}$.

When $j=7$, $G_7\in \mathcal {B}(n-d+s+1,s+1)$. If $x_{i-1}$ lies on
all diametrical paths of $G_7$, then by Lemma \ref{B_{n,d+1}},
$G_7\succ B_{n-d+s+1,s+2}\succ B_{n-d+s+1,s+1}$. Otherwise, in the
same way as in Case 2, we can also obtain $G_7\succ
B_{n-d+s+1,s+1}$.

Similarly, when $j=6$, we can have $G_6\succ B_{n-d+s+1,s+1}$.

When $j=8$, $G_8\in \mathcal {B}(n-d+s,s)$. Since $G_8$ contains no
pendent vertices, then by Lemmas \ref{Bdd_0} and \ref{GB_{n,d+1}},
$G_8\succ B_{n-d+s,s+1}\succ B_{n-d+s,s}$.

{\bf{Subcase 3.2}} There is only one pendent vertex in $G$, say
$x_0$.

Since there are two vertex-disjoint cycles in $G$, $deg_G(x_d)=2$.
Suppose that $x_l$ is the vertex such that $deg_G(x_l)\geq 3$ and
$deg_G(x_i)=2$ for $l+1\leq i\leq d$. It is easy to check that
$l\leq d-2$. Then $G-x_{d-1}x_d\in \mathcal {U}(n,d_9)$, where
$d_9\geq d$, and $G-x_{d-1}-x_d\in \mathcal {U}(n-2,d_{10})$, where
$d_{10}\geq d-1$. Hence,
\begin{eqnarray*}
m_k(G) &=& m_k(G-x_{d-1}x_d)+m_{k-1}(G-x_{d-1}-x_d)\\
       &\geq& m_k(U_{n,d_9})+m_{k-1}(U_{n-2,d_{10}})\\
       &\geq& m_k(U_{n,d})+m_{k-1}(U_{n-2,d-1})\\
       &\geq& m_k(U_{n,d})+m_{k-1}(T_{n-2,d-1})\\
       &\geq& m_k(U_{n,d})+m_{k-1}(P_{d-3}\cup S_{n-d+1})\\
       &=& m_k(B_{n,d}).
\end{eqnarray*}
In particular, $m_2(G)>m_2(B_{n,d})$. Thus, $G\succ B_{n,d}$.

Therefore, we complete the proof. \qed

\begin{theorem}\label{thm4}
Let $G\in \mathcal {B}(n,d)$ with $n\geq 8$, $3\leq d\leq n-3$ and
$G\neq B_{n,d}$. If there is no vertex-disjoint cycles in $G$, then
$G\succ B_{n,d}$.
\end{theorem}
\pf We will prove this theorem by induction on $n-d$.

By Lemma \ref{lem $n-d=3$}, the result holds for $n-d=3$. Let $h\geq
4$ and suppose that the result holds for $n-d<h$. Now assume that
$n-d=h$, let $G\in \mathcal {B}(n,d)$ and $G\neq B_{n,d}$.

{\bf{Case 1}} If $G$ contains no pendent vertices.

Then by Lemmas \ref{Bdd_0} and \ref{GB_{n,d+1}}, $G\succ B_{n,d+1}
\succ B_{n,d}$.

{\bf{Case 2}} If there exists a pendent vertex outside some
diametrical path $P(G)=x_0x_1\ldots x_d$.

Let $u$ be a pendent vertex outside $P(G)$ and $v$ be its unique
neighbor. Then $G-u\in \mathcal {B}(n-1,d)$. If $G-u=B_{n-1,d}$,
then it can be checked that $G-u-v\succ T_{d+1,d-2}$. And thus from
$(*)$, we can obtain that $G\succ B_{n,d}$. Otherwise, by the
induction hypothesis, we have $G-u\succ B_{n-1,d}$. Let $H=G-u-v$,
in order to prove the result, we only need to show that $m_k(H)\geq
m_k(T_{d+1,d-2})$. We choose $C_a$, $C_b$ as above in $G$, and if
there exists the third cycle, denote it by $C_c$.

{\bf{Subcase 2.1}} When $v$ lies on some cycle, say $C_a$.

First, suppose that $v=u_0$ or $u_t$, then $H$ contains no cycles.
If $v$ lies outside $P(G)$, then $H\supseteq P(G)$. Thus $m_k(H)\geq
m_k(P(G))=m_k(P_{d+1})\geq m_k(T_{d+1,d-2})$. If $v$ lies on $P(G)$,
say $v=x_i$.

$(1)$ If $C_a$ and $C_b$ have exactly one common vertex, then
$H\supseteq P_2\cup P_2\cup P_i \cup P_{d-i}$, $P_2\cup P_i\cup
P_{d-i+1}$, $P_2\cup P_{i+1}\cup P_{d-i}$, $P_{i+1}\cup P_{d-i+1}$,
$P_2\cup P_i\cup T_1$, $P_2\cup P_{d-i}\cup T_2$, $P_{i+1}\cup T_1$,
$P_{d-i+1}\cup T_2$ or $T_1\cup T_2$, where $T_1\in \mathcal
{T}(d-i+1,d-i-1)$ and $T_2\in \mathcal {T}(i+1,i-1)$.

If $H\supseteq P_2\cup P_2\cup P_i\cup P_{d-i}$, $P_2\cup P_i\cup
P_{d-i+1}$, $P_2\cup P_{i+1}\cup P_{d-i}$ or $P_{i+1}\cup
P_{d-i+1}$, then $m_k(H)\geq m_k(P_{d+1})\geq m_k(T_{d+1,d-2})$;

If $H\supseteq P_2\cup P_i\cup T_1$, then $m_k(H)\geq
m_k(P_{i+1}\cup T_1)\geq m_k(P_{i+1}\cup T_{d-i+1,d-i-1})\geq
m_k(T_{d+1,d-1})\geq m_k(T_{d+1,d-2})$; Similarly, if $H\supseteq
P_2\cup P_{d-i}\cup T_2$, we also have $m_k(H)\geq
m_k(T_{d+1,d-2})$;

If $H\supseteq T_1\cup T_2$, then $m_k(H)\geq m_k(T_1\cup T_2)\geq
m_k(T_{d-i+1,d-i-1}\cup T_{i+1,i-1})\geq m_k(T_{d+1,d-2})$;

If $H\supseteq P_{i+1}\cup T_1$ or $H\supseteq P_{d-i+1}\cup T_2$,
then $m_k(H)\geq m_k(T_1\cup T_2)\geq m_k(T_{d+1,d-2})$.

$(2)$ If $C_a$ and $C_b$ have at least two common vertices, then
$H\supseteq P_3\cup P_i\cup P_{d-i}$, $P_i\cup P_{d-i+2}$, $P(G)$,
$P_i\cup T_3$ or $P_i\cup T_4$, where $T_3\in \mathcal
{T}(d-i+2,d-i-1)$ and $T_4\in \mathcal {T}(d-i+2,d-i)$.

If $H\supseteq P_3\cup P_i\cup P_{d-i}$, $P_i\cup P_{d-i+2}$ or
$P(G)$, then $m_k(H)\geq m_k(P_{d+1})\geq m_k(T_{d+1,d-2})$;

If $H\supseteq P_i\cup T_3$, then $m_k(H)\geq m_k(P_i\cup T_3)\geq
m_k(P_i\cup T_{d-i+2,d-i-1}) \geq m_k(T_{d+1,d-2})$;

If $H\supseteq P_i\cup T_4$, then $m_k(H)\geq m_k(P_i\cup T_4)\geq
m_k(P_i\cup T_{d-i+2,d-i})\geq m_k(T_{d+1,d-1})\geq
m_k(T_{d+1,d-2})$.

Next, suppose that $v\neq u_0$ and $v\neq u_t$. If $v$ lies outside
$P(G)$, then $H\supseteq P(G)$, similarly, we have $m_k(H)\geq
m_k(T_{d+1,d-2})$. So suppose that $v$ lies on $P(G)$. Then $P(G)$
and $C_a$ have common vertices, say $x_i,\ldots,x_{i+p}$, where
$p\geq 0$.

$(1)$ If $p=0$, then $i\geq 1$ and $H\supseteq P_i\cup P_{d-i}\cup
C_s$, where $s=b$ or $c$. It follows that
\begin{eqnarray*}
m_k(H)&\geq& m_k(P_i\cup P_{d-i}\cup C_s)\geq m_k(P_i\cup
P_{d-i}\cup S_s)\geq
m_k(P_i\cup P_{d-i}\cup P_3)\\
&\geq& m_k(P_{d+1})\geq m_k(T_{d+1,d-2}).
\end{eqnarray*}

$(2)$ If $p\geq 1$. When $v\neq x_i$ and $v\neq x_{i+p}$, then
$H\supseteq G_1$, where $G_1\in \mathcal {U}(s_1,d_1)$, $d_1\geq d$
and $s_1\geq d_1+2$. Thus
$$m_k(H)\geq m_k(G_1)\geq m_k(U_
{s_1,d_1})\geq m_k(T_{s_1,d_1})\geq m_k(T_{s_1,d})\geq
m_k(T_{d+1,d-2}).$$ When $v=x_i$ or $x_{i+p}$, say $v=x_i(i\geq 1)$,
then $H\supseteq P_i\cup G_2$, where $G_2\in \mathcal {U}(s_2,d_2)$,
$d_2\geq
 d-i-1\geq 2$ ($d-i-1\geq 1$, for $d-i-1=1$, clearly $m_k(H)\geq m_k(T_{i+3,i})$) and $s_2\geq d-i+2$;
or $H\supseteq P_i\cup G_3$, where $G_3\in \mathcal {U}(s_3,d_3)$,
$d_3\geq d-i\geq 2$ and $s_3\geq d_3+2$; or $H\supseteq P_i\cup G'$,
where $G'$ is the graph obtained by attaching a path $P_{d-i-2}$ to
a vertex of $C_b=C_3$.

Suppose that $H\supseteq P_i\cup G_2$. If $d_2=2$, then $i=d-3$ and
$s_2\geq 5$. Hence,
$$m_k(H)\geq
m_k(P_i\cup G_2)\geq m_k(P_i\cup S_{s_2})\geq m_k(P_i\cup
T_{5,2})\geq m_k(T_{i+4,i+1})=m_k (T_{d+1,d-2}).$$ If $d_2\geq 3$,
then
\begin{eqnarray*}
m_k(H)&\geq& m_k(P_i\cup G_2)\geq m_k(P_i\cup U_{s_2,d_2})\geq m_k(P_i\cup T_{s_2,d_2})\\
&\geq& m_k(T_{i+s_2-1,i+d_2-1})\geq m_k(T_{d+1,d-1})\geq
m_k(T_{d+1,d-2}).
\end{eqnarray*}

Suppose that $H\supseteq P_i\cup G_3$. If $d_3=2$, then $i=d-2$ and
$s_3\geq 4$. Thus
\begin{eqnarray*}
m_k(H)&\geq& m_k(P_i\cup G_3)\geq m_k(P_i\cup S_{s_3})\geq m_k(P_i\cup T_{4,2})\\
&\geq& m_k(T_{d+1,d-1})\geq m_k(T_{d+1,d-2}).
\end{eqnarray*}
If $d_3\geq 3$, then
\begin{eqnarray*}
m_k(H)&\geq& m_k(P_i\cup G_3)\geq m_k(P_i\cup U_{s_3,d_3})\geq m_k(P_i\cup T_{s_3,d_3})\\
&\geq& m_k(T_{i+s_3-1,i+d_3-1})\geq m_k(T_{d+1,d-1})\geq
m_k(T_{d+1,d-2}).
\end{eqnarray*}

Suppose that $H\supseteq P_i\cup G'$. If $d-i-2=0$, then $i=d-2$ and
$H\supseteq P_i\cup C_3$. Therefore,
$$m_k(H)\geq m_k(P_i\cup G')\geq m_k(P_i\cup T_{4,2})\geq m_k(T_{i+3,i+1})=m_k
(T_{d+1,d-1})\geq m_k(T_{d+1,d-2}).$$ If $d-i-2\geq 1$, choose an
edge $u_0u_1$ of $G'$ in $C_3$ such  that both $u_0$ and $u_1$ have
degree $2$. Then
\begin{eqnarray*}
m_k(H)&\geq& m_k(P_i\cup G')\\
&=&m_k(P_i\cup G'-u_0u_1)+
m_{k-1}(P_i\cup G'-u_0-u_1)\\
&=&m_k(P_i\cup T_{d-i+1,d-i-1})+m_{k-1}(P_i\cup P_{d-i-1})\\
&\geq& m_k(T_{d,d-2})+m_{k-1}(P_{d-3})\\
&=&m_k(T_{d+1,d-2}).
\end{eqnarray*}

{\bf{Subcase 2.2}} When $v$ lies outside any cycle.

In this case, $H$ contains two cycles $C_a$ and $C_b$ with at least
one common vertex. Let $C_a\cdot C_b$ denote the subgraph of $G$
induced by $V(C_a)\cup V(C_b)$.

First, suppose that $v$ lies on $P(G)$, say $v=x_i$.

If vertices on $P(G)$ lie outside any cycle, then $H\supseteq
C_a\cdot C_b\cup P_i\cup P_{d-i}$. Thus,
$$m_k(H)\geq m_k(C_a\cdot C_b\cup P_i\cup
P_{d-i})\geq m_k(P_3\cup P_i\cup P_{d-i})\geq m_k(P_{d+1})\geq
m_k(T_{d+1,d-2}).$$ Otherwise, $H\supseteq P_i\cup G_1$, where
$G_1\in \mathcal {B}(s_1,d_1)$, $d_1\geq \max\{d-i-1,2\}$ and
$d_1+2\leq s_1\leq n-2-i$.

Suppose that $s_1\geq d_1+3$. If $d_1=2$, then $i\geq d-3$ and
$s_1\geq 5$. Thus we have
$$m_k(H)\geq m_k(P_i\cup G_1)\geq m_k(P_i\cup S_{s_1})\geq m_k(P_i\cup T_{5,2})\geq m_k(T_{i+4,i+1})\geq m_k(T_{d+1,d-2}).$$
If $d_1\geq 3$, then $d-i-1\leq d_1$, which deduces that
$s_1-d_1\leq n-2-i- (d-i-1)=n-d-1<h$. By the induction hypothesis,
$G_1\succ B_{s_1,d_1}$, therefore,
$$m_k(H)\geq m_k(P_i\cup G_1)\geq m_k(P_i\cup B_{s_1,d_1})\geq m_k(P_i\cup T_{s_1,d_1})\geq m_k(T_{d+1,d-2}).$$

Now suppose that $s_1=d_1+2$. In this case, $G_1$ is obtained by
attaching respectively paths $P_l(0\leq l\leq d_1-2)$ and
$P_{d_1-l-2}$ to the two non-adjacent vertices in $K_4-e$. If
$d_1=2$, then $i\geq d-3$. It can be easily checked that
$m_k(K_4-e)\geq m_k(T_{5,2})$. Thus $m_k(H)\geq m_k(P_i\cup
(K_4-e))\geq m_k(P_i\cup T_{5,2})\geq m_k(T_{d+1,d-2})$. If $d_1\geq
3$, choose an edge $u_0u_1$ of $G_1$ such that $u_0$ and $u_1$ are
both of degree $3$ in $K_4-e$. Then we get
\begin{eqnarray*}
m_k(H)&\geq& m_k(P_i\cup G_1)\\
&=&m_k(P_i\cup G_1-u_0u_1)+m_{k-1}(P_i\cup
G_1-u_0-u_1)\\
&\geq& m_k(P_i\cup U_{s_1,d_1})+m_{k-1}(P_i\cup P_{l+1}\cup P_{d_1-l-1})\\
&\geq& m_k(P_i\cup T_{s_1,d_1})+m_{k-1}(P_i\cup P_{d_1-1})\\
&\geq& m_k(P_i\cup T_{d-i+1,d-i-1})+m_{k-1}(P_i\cup P_{d-i-2})\\
&\geq& m_k(T_{d,d-2})+m_{k-1}(P_{d-3})\\
&=&m_k(T_{d+1,d-2}).
\end{eqnarray*}

Next, suppose that $v$ lies outside $P(G)$.

In this case, $H\supseteq G_2$ or $H\supseteq C_a\cdot C_b\cup
P(G)$, where $G_2\in \mathcal {B}(s,d)$ with $d+2\leq s\leq n-2$. It
is easy to verify that $m_k(H)\geq m_k(T_{d+1,d-2})$.

{\bf{Case 3}} If any diametrical path of $G$ contains all pendent
vertices in $G$.

We can obtain that $G\succ B_{n,d}$ by similar arguments as those in
Case 3 of Theorem 2.

Consequently, the proof is complete. \qed

Combining Theorems \ref{thm3} and \ref{thm4}, we obtain the
following main result of this section.

\begin{theorem}\label{thm5}
Let $G\in \mathcal {B}(n,d)$ with $n\geq 8$, $3\leq d\leq n-3$ and
$G\neq B_{n,d}$. Then $ME(G)> ME(B_{n,d})$.
\end{theorem}
\pf According to Theorems \ref{thm3} and \ref{thm4}, we have known
that $G\succ B_{n,d}$. And then using the increasing property
(namely, $G_1\succ G_2\Longrightarrow ME(G_1)> ME(G_2)$), we get the
result we want. \qed

\begin{figure}[h,t,b,p]
\begin{center}
\includegraphics[scale = 1.0]{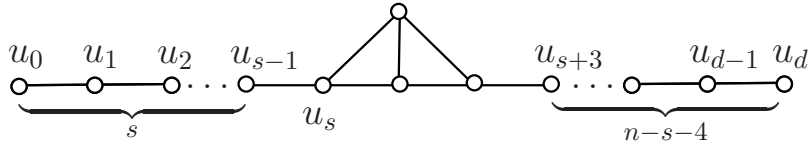}
\caption{Graph $B_n^s$ with $0\leq s \leq \lfloor
n/2\rfloor-2$.}\label{fig d=n-2}
\end{center}
\end{figure}

We will conclude this section by discussing the case $d=n-2$. Since
any graph $G$ in $\mathcal {B}(n,n-2)$ is of the form $B_n^s$(as
shown in Figure \ref{fig d=n-2}), where $0\leq s \leq \lfloor
n/2\rfloor-2$. Through simple analysis, we get the following result.

\begin{theorem}\label{thm6}
Let $G\in \mathcal {B}(n,n-2)$ with $n\geq 6$ and $G\neq B_n^1$,
then $ME(G)> ME(B_n^1)$.
\end{theorem}
\pf Since $G\in \mathcal {B}(n,n-2)$ and $G\neq B_n^1$, then $G$ is
$B_n^s$, where $s=0,2,\ldots,\lfloor n/2\rfloor-2$.

{\bf{Case 1}} $G=B_n^0$.

Then $G-u_d=B_n^0-u_d=B_n^1-u_0$ and
$G-u_{d-1}-u_d=B_n^0-u_{d-1}-u_d\succ B_n^1-u_0-u_1$. Thus we have
$G=B_n^0\succ B_n^1$ by Lemma \ref{lem9}.

{\bf{Case 2}} $G=B_n^s(s=2,\ldots,\lfloor n/2\rfloor-2)$.

Then
\begin{eqnarray*}
m_k(G=B_n^s)&=& m_k(B_n^s-u_{s-2}u_{s-1})+m_{k-1}(B_n^s-u_{s-2}-u_{s-1})\\
            &=& m_k(P_{s-1}\cup B_{n-s+1}^1)+m_{k-1}(P_{s-2}\cup B_{n-s}^0),
\end{eqnarray*}
together with
\begin{eqnarray*}
m_k(B_n^1)&=& m_k(B_n^1-u_{d-s+1}u_{d-s+2})+m_{k-1}(B_n^1-u_{d-s+1}-u_{d-s+2})\\
          &=& m_k(P_{s-1}\cup B_{n-s+1}^1)+m_{k-1}(P_{s-2}\cup B_{n-s}^1).
\end{eqnarray*}

By Case $1$, we have got that $B_{n-s}^0\succ B_{n-s}^1$, then
$m_{k-1} (P_{s-2}\cup B_{n-s}^0)\geq m_{k-1}(P_{s-2}\cup B_{n-s}^1)$
and $m_{2}(P_{s-2}\cup B_{n-s}^0)> m_{2}(P_{s-2}\cup B_{n-s}^1)$.
Thus $G=B_n^s\succ B_n^1$.

Therefore, we always have $G\succ B_n^1$. And then $ME(G)>
ME(B_n^1)$. \qed

\section{Summary}
In \cite{ISG2,SG1}, the authors introduced the concept of
``set-complexity", based on a context-dependent measure of
information, and used this concept to describe the complexity of
gene interaction networks. The binary graphs and edge-colored graphs
are studied and the relation between complexity and structure of
these graphs is examined in detail. In contrast, we put the emphasis
on analyzing properties of spectra-based entropies and study
interrelations thereof.

In this paper, we characterize the graphs with minimal matching
energy among all unicyclic and bicyclic graphs with a given diameter
$d$. With respect to matching energy of graphs, $U_{n,d}$ and
$B_{n,d}$ are two extremal graphs in $\mathcal {U}(n,d)$ and
$\mathcal {B}(n,d)$ respectively. Moreover, both of them are
interesting and have the similar extremum property in other aspects.
For example, among all unicyclic graphs of a given diameter,
$U_{n,d}$ is the extremal graph on graph energy \cite{LZ}. Besides,
it is also the underling graph of the extremal graph on skew energy
\cite{YGX}. In addition, $B_{n,d}$ has the minimal energy in one
class of bicyclic graphs with a given diameter \cite{YZ}. From this
point, we guess that this two graphs may also be the extremal graphs
on some other parameters of graphs. Studying the properties of this
two graphs will be one of the future work of us. An important
question is how general the bounds are. Obviously, the proof
techniques use structural properties of the graphs under
consideration and it may be intricate to extend the techniques when
using more general graphs. On the other hand, the roots of graph
polynomials could be used to characterize graphs structurally.
This will be one of the future work. For more results, we refer to \cite{DEG, KDE}.\\

\noindent{\bf Acknowledgement.} The authors would like to thank the
referees for valuable comments. The authors are supported by NSFC,
PCSIRT, China Postdoctoral Science Foundation (2014M551015) and
China Scholarship Council.

\end{document}